# AN ESTIMATE FOR THE MULTIPLICITY OF BINARY RECURRENCES

ROBERTO GIORGIO FERRETTI

ABSTRACT

*We use a refined version of Roth's Lemma, proved with the help of Faltings' Product Theorem, in order to give an upper bound for the multiplicity of a binary linear recurrence.*

Résumé

*Nous utilisons une version du Lemme de Roth, provenant du Théorème du Produit de Faltings, pour majorer la multiplicité d'une récurrence linéaire binaire.*

Version française abrégée

Dans l'article [7] H.-P. Schlickewei donne une majoration, indépendante du corps des nombres, de la multiplicité d'une récurrence linéaire binaire. La contribution principale à cette majoration est exprimée par une version, pour la droite projective, du théorème du sous-espace de Schmidt. Nous améliorons cette majoration en utilisant la même méthode, mais faisant appel à un Lemme de Roth plus puissant, démontré à partir d'idées liées au théorème du Produit de Faltings ([4], [3], [6]). Même si ce résultat n'est pas comparable à la majoration remarquable obtenue récemment par F. Beukers et H.-P. Schlickewei [2], elle montre jusqu'à quel point on peut arriver en utilisant les techniques du théorème du sous-espace de Schmidt dans ce contexte. Nous démontrons,

**Théorème 1.1.** *Soient $a, b, \alpha, \beta$ nombres complexes, tels que au moins un entre $\alpha, \beta$ n'est pas une racine de l'unité. Alors il y a au moins*

$$2^{57}$$

*entiers $m \in \mathbb{Z}$ tels que*

$$a\alpha^m + b\beta^m + 1 = 0.$$

Pour démontrer ce théorème il nous faut la proposition suivante. Soit $S$ un sous-ensemble fini de $M_K$ contenant les places à l'infini. Pour tout $v \in S$ nous donnons deux formes linéaires $L_{1,v}(x), L_{2,v}(x) \in \{X_1, X_2, X_1 + X_2\}$ et deux nombres réels $e_{1,v}, e_{2,v}$ tels que $\sum_{v \in S}(e_{1,v} + e_{2,v}) = 0$, et pour tout sous-ensembles $S'$ of $S$, et chaque $(i(v))_{v \in S'}$ avec $i(v) \in \{1, 2\}$, $|\sum_{v \in S'} e_{i(v),v}| \leqslant 1$. Ce sont les conditions $(3.1) - (3.3)$ de [7]. On obtient,

**Proposition 2.1.** *Supposons que $L_{v,i}, e_{v,i}$ soient comme en haut. Soit $0 < \delta < 1$ un nombre réel. Alors pour tout nombres réels $Q > 4^{\delta}$ les solutions $x \in K^2$ de*

$$\begin{aligned}
\|L_{i,v}(x)\|_v &< Q^{e_{i,v} - \delta[K_v:K]/[K:\mathbb{Q}]}, \quad v \in S,\ v|\infty,\ i \in \{1,2\}, \\
\|L_{i,v}(x)\|_v &\leqslant Q^{e_{i,v}}, \quad v \in S,\ v \nmid \infty,\ i \in \{1,2\}, \\
\|x\|_v &\leqslant 1, \quad v \in S,
\end{aligned}$$

*sont contenues dans au plus $2^{227/10} \delta^{-3} \log \delta^{-1}$ droites de $K^2$.*

Nous rappelons que l'indice $i_{x,\mathbf{r}}(P)$ d'un polynôme $P$ en $2m$ variables, multihomogène de multidegré $\mathbf{r} = (r_1, \cdots, r_m) \in \mathbb{Z}^m$ en un point $x \in K^{2m}$ est la multiplicité de $P$ en $x$ calculée avec le poids $1/r_i$ dans la $i$-ème direction. Nous avons utilisé ici la version suivante du Lemme de Roth.

**Lemme 3.1. (Lemme de Roth)** *Soit $m \geqslant 2$ un entier, $\mathbf{r} = (r_1, \cdots, r_m) \in Zz^m$ des entiers positifs et $0 < \vartheta \leqslant m^2(m+1)$ un nombre réel. Supposons que pour tout entier $h$ avec $1 \leqslant h \leqslant m-1$, $r_h/r_{h+1} \geqslant \frac{m^2(m+1)}{\vartheta}$. Soit $P$ un polynôme non-nul en $2m$ variables, multi-homogène de multidegré $\mathbf{r}$ et soient $L_1, \cdots, L_m$ des formes lineaires binaires avec coefficients en $\mathbb{Q}$ tels que, pour tout entiers $h$ avec $1 \leqslant h \leqslant m$*

$$r_h \log H(L_h) \geqslant \frac{7m(m!)^2 m^m}{2\vartheta^m} (\sum_{i=1}^{m} r_i + \log H(P)),$$

*Alors, pour chaque entier $h$ avec $1 \leqslant h \leqslant m$, il y a un $x_h \in V(L_h)(K)$ tel que pour $x = (x_1, \cdots, x_m)$ nous avons $i_{x,\mathbf{r}}(P) < \vartheta$.*

Le théorème principal decoule maintenant comme dans [7]. Nous utilisons le même principe des trous et les mêmes choix des paramètres.

## 1. Introduction

In his paper [7] H.-P. Schlickewei gives an upper bound, independent from the degree of the number field, for the multiplicity of a binary linear recurrence. The main contribution to his bound comes from a particular version of the Schmidt Subspace Theorem in the case of the projective line. We improve the upper bound for the multiplicity of a binary linear recurrence using the same method but applying a refined version of Roth's Lemma, which combines the ideas of [3] and [6] arising from the Faltings' Product Theorem ([4]). Although this result falls short of the remarkable bound recently obtained by F. Beukers and H.-P. Schlickewei [2], this method still gives better bounds compared to the result of E. Bombieri, J. Müller and M. Poe [1]. We obtain,

**Theorem 1.1.** *Let $a, b, \alpha, \beta$ be complex numbers, such that one at least among $\alpha, \beta$ is not a root of unity. Then there are at most*

$$2^{57}$$

*integers $m \in \mathbb{Z}$ such that*

(1) $$a\alpha^m + b\beta^m + 1 = 0.$$

**Corollary 1.2.**   1. *For $0 \neq \nu_0, \nu_1$ complex numbers let*

$$u_{n+2} = \nu_1 u_{n+1} + \nu_0 u_n$$

*be a binary linear recurrence sequence of complex numbers. Suppose that one at least of the roots $\alpha_1$, $\alpha_2$ of the polynomial $z^2 - \nu_1 z - \nu_0$ is not a root of unity, and if $\alpha_1 \neq \alpha_2$ assume moreover that $\alpha_1/\alpha_2$ is not a root of unity. Then $\{u_n\}$ has multiplicity*

$$U \leqslant 2^{57}.$$

2. *Let $0 \neq \mu_0$, $\mu_1$, $\mu_2$ complex numbers, let us consider a ternary sequence*

$$v_{m+3} = \mu_2 v_{m+2} + \mu_1 v_{m+1} + \mu_0$$

*with $|v_0| + |v_1| + |v_2| \neq 0$. If the polynomial $z^3 - \mu_2 z^2 - \mu_1 z - \mu_0$ has three distinct roots $\alpha_1$, $\alpha_2$, $\alpha_3$, assume that at least one of the quotients $\alpha_1/\alpha_3$, $\alpha_2/\alpha_3$ is not a root of unity. Then $\{v_n\}$ has zero-multiplicity*

$$U(0) \leqslant 2^{57}.$$

*Proof.* We prove this corollary as in [7] §1 using Theorem 1.1 instead of [7] Theorem 1.   □

Here for a complex number $c$ we define the *c-multiplicity* $U(c)$ of $\{u_n\}$ as the number of solutions $m \in \mathbb{Z}$ of the equation $u_m = c$ and we write $U = \sup_c U(c)$ and call $U$ the *multiplicity* of the sequence.

Theorem 1 of [7] differs from Theorem 1.1 in that instead of $2^{57}$ it has

$$2^{2^{23}}$$

as upper bound for the number of solutions of (1). We refer to [7] for a discussion concerning known results and conjectures on linear recurrences.

## 2. Subspace Theorem

Let $K$ be an algebraic number field. Denote its ring of integers by $\mathcal{O}_K$ and its collection of places (equivalence classes of absolute values) by $M_K$. For $v \in M_K$, $x \in K$ we define the absolute value $|x|_v$ by

1. $|x|_v = |\sigma(x)|^{1/[K:\mathbb{Q}]}$ if $v$ corresponds to a real embedding $\sigma : K \hookrightarrow \mathbb{C}$,
2. $|x|_v = |\sigma(x)|^{2/[K:\mathbb{Q}]}$ if $v$ corresponds to a pair of conjugate complex embeddings $\sigma, \overline{\sigma} : K \hookrightarrow \mathbb{C}$,
3. $|x|_v = (N\mathcal{P})^{-ord_\mathcal{P}(x)/[K:\mathbb{Q}]}$ if $v$ corresponds to the prime ideal $\mathcal{P}$ of $\mathcal{O}_K$.

Here $N\mathcal{P} = \#(\mathcal{O}_K/\mathcal{P})$ is the norm of $\mathcal{P}$ and $ord_\mathcal{P}(x)$ is the exponent of $\mathcal{P}$ in the prime ideal decomposition of the principal ideal generated by $x$, the order of $0$ is $\infty$. We denote by $K_v$ the algebraic closure of the $v$-adic completion of $K$. In the first two cases we call $v$ infinite and write $v \mid \infty$, in case 3 we call $v$ finite and write $v \nmid \infty$. These absolute values satisfy the product formula $\prod_v |x|_{v \in M_K} = 1$ for $x \in K^*$. If $x = (a_1, \cdots, a_m) \in K^m \setminus \{0\}$ we put

$$\begin{aligned}
\|x\|_v &= (\sum_{i=1}^m |a_i|_v^{2[K:\mathbb{Q}]})^{1/2[K:\mathbb{Q}]}, &&\text{if } v \text{ is real,} \\
\|x\|_v &= (\sum_{i=1}^m |a_i|_v^{[K:\mathbb{Q}]})^{1/[K:\mathbb{Q}]}, &&\text{if } v \text{ is complex,} \\
\|x\|_v &= \max\{|a_1|, \cdots, |a_m|\}, &&\text{if } v \nmid \infty.
\end{aligned}$$

Now define the *height* of $x$ as
$$H(x) = \prod_{v \in M_K} \|x\|_v.$$

By the product formula this defines a function on the projective space $\mathbb{P}^{m-1}(K)$. Further it depends only on the point $x$ and not on the choice of the number field $K$ containing the coordinates of $x$. For a linear form $L(x) = a_1 x_1 + \cdots + a_n x_n$ with algebraic coefficients, we define the height $H(L)$ as the height of the point $(a_1, \cdots, a_m)$. Moreover we define the height $H(V(L))$ of the $(n-1)$-dimensional linear subspace of $K^n$
$$V(L) = \{x \in K^n : L(x) = 0\}$$
as the height of the linear form $L$. Similarly the height of a polynomial is the height of the sequence of its coefficients.

Let us consider for $n = 2$ the set of linear forms given by
$$\mathbf{L} = \{L_1(x) = x_1, \quad L_2(x) = x_2, \quad L_3(x) = x_1 + x_2\}.$$

Let $S$ be a finite subset of $M_K$ containing all infinite places. We suppose that for each $v \in S$ we are given a pair of different linear forms $L_{1,v}(x), L_{2,v}(x)$ out of $\mathbf{L}$ and a pair of real numbers $e_{1,v}, e_{2,v}$ such that

(2) $$\sum_{v \in S}(e_{1,v} + e_{2,v}) = 0,$$

and for each subset $S'$ of $S$, and any tuple $(i(v))_{v \in S'}$ with $i(v) \in \{1, 2\}$,

(3) $$|\sum_{v \in S'} e_{i(v),v}| \leqslant 1.$$

These are the conditions $(3.1) - (3.3)$ of [7]. Consider, for a real number $0 < \delta < 1$, the simultaneous inequalities

(4) $$\begin{aligned}\|L_{i,v}(x)\|_v &< Q^{e_{i,v} - \delta[K_v:K]/[K:\mathbb{Q}]}, & v \in S, \ v|\infty, \ i \in \{1,2\}, \\ \|L_{i,v}(x)\|_v &\leqslant Q^{e_{i,v}}, & v \in S, \ v \nmid \infty, \ i \in \{1,2\}, \\ \|x\|_v &\leqslant 1, & v \in S,\end{aligned}$$

which correspond to condition (3.4) of [7]. We are now ready to state our improvement of Schlickewei's Lemma 3.1 [7].

**Proposition 2.1.** *Suppose that $\delta, L_{v,i}, e_{v,i}$ are as above and that (2), (3) hold. Then for all real numbers*
$$Q > 4^\delta$$
*the solutions $x \in K^2$ of (4) are contained in the union of at most*
$$2^{227/10} \delta^{-3} \log \delta^{-1}$$
*one-dimensional linear subspaces of $K^2$.*

## 3. Roth's Lemma

Let $P(x_{11}, x_{12}; \cdots ; x_{m1}, x_{m2})$ be a polynomial with rational coefficients, multihomogeneous of multidegree $\mathbf{r} = (r_1, \cdots, r_m)$. Given a point $x \in K^m$ we define the *index* $i_{x,\mathbf{r}}(P)$ of $P$ with respect to $(x, \mathbf{r})$ as the weighted multiplicity of $P$ at $x$ with weights $1/d_j$. This turns out to be the same as defining the index with respect to binary linear forms as in [7] §5.

**Lemma 3.1. (Roth's Lemma)** *Let $m \geqslant 2$ be an integer, $\mathbf{r} = (r_1, \cdots, r_m)$ a $m$-tuple of positive integers and $0 < \vartheta \leqslant m^2(m+1)$ a real number. Suppose that for all integers $h$ with $1 \leqslant h \leqslant m-1$*

$$(5) \qquad r_h/r_{h+1} \geqslant \frac{m^2(m+1)}{\vartheta}.$$

*Furthermore, let $P$ be a non-zero polynomial in $2m$ variables, multihomogeneous of multidegree $\mathbf{r}$ and let $L_1, \cdots, L_m$ be binary linear forms with coefficients in $\mathbb{Q}$ such that, for all integers $h$ with $1 \leqslant h \leqslant m$*

$$(6) \qquad r_h \log H(L_h) \geqslant \frac{7m(m!)^2 m^m}{2\vartheta^m}(\sum_{i=1}^{m} r_i + \log H(P)),$$

*Then, for all integers $h$ with $1 \leqslant h \leqslant m$, there is a $x_h \in V(L_h)(K)$ such that for $x = (x_1, \cdots, x_m)$ we have*

$$i_{x,\mathbf{r}}(P) < \vartheta.$$

*Proof.* We follow the proof of [3] Theorem 3 avoiding [7] Lemma 9 but using [6] Proposition 5.3 instead. Let $h$ be an integer with $1 \leqslant h \leqslant m-1$ then

$$r_h \log H(V(L_h)) \leqslant \frac{m!(m+1)^m}{\vartheta^m}(\frac{m!}{2}\sum_{i=1}^{m} r_i + m!(\log H(P) + \log(2^{-1+\sum_{i=1}^{m} r_i}(\sum_{i=1}^{m} r_i)^{m-1})))$$

The difference comes out in [3] (5.3) where we obtain the inequality

$$r_h \log H(V(L_h)) \leqslant \frac{(m!)^2 m^m e}{\vartheta^m}((\frac{1}{2} + \log 2)\sum_{i=1}^{m} r_i + (m-1)\log(\sum_{i=1}^{m} r_i) + \log H(P)) \leqslant$$
$$\leqslant \frac{7m(m!)^2 m^m}{2\vartheta^m}(\sum_{i=1}^{m} r_i + \log H(P)).$$

Since $H(L_h) = H(V(L_h))$ this completes the proof of Lemma 3.1. $\square$

## 4. Proof of the Subspace Theorem

Let $\mathbf{g}_1(Q)$, $\mathbf{g}_2(Q)$, $V_h(Q)$, for $h$ integer with $1 \leqslant h \leqslant m$, be as in [7] §5. Then we have the following version of [7] Lemma 6.1,

**Lemma 4.1.** *Suppose $0 < \delta < 1$ is real and that*

$$m > 28800\delta^{-2}.$$

Put

$$E = \frac{m^2(m+1)}{240}, \quad F = \frac{7}{2}m(m!)^2(\frac{m}{480})^m.$$

*Then the numbers $Q$ with*

$$L_1(\mathbf{g}_1(Q))L_2(\mathbf{g}_1(Q))L_3(\mathbf{g}_1(Q)) \neq 0, \tag{7}$$

$$\lambda_1(Q) < Q^{-\delta}, \tag{8}$$

$$Q^{\delta^2} > 2^{600mF}, \tag{9}$$

*are contained in the union of at most $m-1$ intervals of the type*

$$Q_h < Q \leqslant Q_h^E.$$

*Proof.* The proof goes along the same lines as the proof of [7] Lemma 6.1. We construct as in loc. cit. the numbers $Q_h$ for $1 \leqslant h \leqslant m$ with $Q_{h+1} \geqslant Q_h^E$, we choose as there $\varepsilon = \delta/60$ and put $r_h = [r_1 \log Q_1/\log Q_k] + 1$. Then the polynomial $P$ of [7] Lemma 5.1 has index $\geqslant m\varepsilon$ with respect to $(V_1(Q_1), \cdots, V_m(Q_m); \mathbf{r})$. Put $\vartheta = 480$. By construction, for all integers $h$ with $1 \leqslant h \leqslant m-1$, we have $r_h/r_{h+1} \geqslant E/2 = m^2(m+1)\vartheta^{-1}$, so (5) is satisfied. If $\Gamma = \delta/10$, then we obtain as in [7] Lemma 6.1 that

$$r_h \log H(V_h) > r_1 \Gamma \log H(V_1).$$

Moreover by [7] Lemma 4.2 and (9),

$$\log H(V_h) \geqslant \Gamma \log H(V_h) > \Gamma^2 \log Q_h = \frac{\delta^2}{100} \log Q_h > 6mF \log 2. \tag{10}$$

By [7] Lemma 5.1 $P$ has $\log H(P) < 4mr_1 \log 2$. Finally by combining this with (10) we get

$$\frac{7m(m!)^2 m^m}{2\vartheta^m}(\sum_{j=1}^m r_i + \log H(P)) < F(mr_1 + 4mr_1 \log 2) < 6mFr_1 \log 2 < \frac{\delta^2 r_1}{100} \log Q_h <$$
$$< \Gamma r_i \log H(V_1) \leqslant r_h \log H(V_h).$$

The conclusion of Lemma 3.1 is that there is a point $x = (x_1, \cdots, x_m) \in K^m$ such that for all integers $h$ with $1 \leqslant h \leqslant m$, $x_h \in V_h$ and

$$\mathrm{ind}_{x,\mathbf{r}}P < \vartheta = 480 < \frac{8.60}{\delta} = \frac{8}{\varepsilon} < m\varepsilon.$$

This yields the desired proof. $\square$

The proof of Proposition 2.1 now follows easily. It is clear from [7] §6 that there are no more than

$$m(1 + \frac{4}{\delta}\log E) + (1 + \frac{4}{\delta}\log(300\delta^{-1}F)) \tag{11}$$

subspaces. Since $0 < \delta < 1$ we can choose $m \leqslant 28801.\delta^{-2}$, then we get

$$\log E < 26 + 6\log \delta^{-1},$$

and by Stirling formula

$$\log(300\delta^{-1}F) < 767865\delta^{-2}\log \delta^{-1}.$$

The bound (11) does not exceed $2^{227/10}\delta^{-3}\log \delta^{-1}$, and this confirms Proposition 2.1.

## 5. Proof of the Main Theorem

We follow here the same proof as [7] §10. We use $\delta = 1/9$. By Proposition 2.1 for values of $Q > 4^9$, given any pair $(i(v), (e_{i,v}))$ of [7] Lemma 8.1, the solutions of (1) are contained in the union of not more than

$$2^{227/10}.3^6 \log 3$$

one-dimensional linear subspaces. The other constants that are relevant for our estimate remain the same as in [7] §10, thus the number of solutions of (1) is bounded above by

$$2(2^7.4800.36 \log 4 + 1 + 2^{227/10}.6^{12} \log 3) < 2^{57}.$$

This proves Theorem 1.1.

Aknowledgements: This work was carried out with support of Swiss National Research Foundation and Fondazione Stefano Franscini. The author benefited from the hospitality of I.H.E.S. and useful discussions with J.-B. Bost, C. Soulé and G. Wüstholz.

INSTITUT DES HAUTES ÉTUDES SCIENTIFIQUES, 35, ROUTE DE CHARTRES, F-91440 BURES-SUR-YVETTE
*E-mail address*: `ferretti@ihes.fr`